\documentclass[12pt]{article}
\usepackage{amssymb, amsmath, url, graphicx}

\def\3{\subset }
\def\4{\subseteq }
\def\<{\left<}
\def\>{\right>}

\def\bit{\begin{itemize}}
\def\eit{\end{itemize}}
\def\3{\subset }
\def\4{\subseteq }

\def\calc{{\cal C}}

\def\0{\leqno}

\def\barr{\begin{array}}
\def\earr{\end{array}}

\def\Z{{\rlap{$\kern2pt{\rm Z}$}{\rm Z}\,}}

\title{A note on a class of gyrogroups}
\author{Marius T\u arn\u auceanu}
\date{October 1, 2016}

\begin{document}

\maketitle

\begin{abstract}
In \cite{1}, a class of gyrogroups $(G,\odot)$ has been associated to certain groups $(G,\cdot)$.
We give a necessary and sufficient condition for $(G,\odot)$ to be gyrocommutative. We also prove
that under a suitable assumption two finite groups central by a $2$-Engel group are isomorphic if and only if their
associated gyrogroups are isomorphic.
\end{abstract}



\section{Introduction}

Gyrogroups are suitable generalization of groups, whose origin is
described in \cite{7,8}. They share remarkable analogies with
groups. In fact, eve\-ry group forms a gyrogroup under the same operation.
Many of classical theorems in group theory also hold for gyrogroups,
including the Lagrange theorem \cite{3}, the fundamental isomorphism theorems
\cite{4}, and the Cayley theorem \cite{4} (for all these theorems see also \cite{6}).
Gyrogroup actions and related results, such as the orbit-stabilizer theorem,
the orbit decomposition theorem, and the Burnside lemma have been studied in \cite{5}.

The present note deals with a connection between groups and gyrogroups, namely with
the gyrogroup associated to any group central by a 2-Engel group (see \cite{1}). We determine 
conditions for such a gyrogroup to be gyrocommutative and for such two gyrogroups to be isomorphic. An interesting
conservative functor between a subcategory of groups and the category of gyrogroups is constructed, as well.

Recall that a groupoid $(G,\odot)$ is called a \textit{gyrogroup} if its binary operation
satisfies the following axioms:
\begin{itemize}
\item[1.] There is an element $e\in G$ such that $e\odot a=a$ for all $a\in G$.
\item[2.] For every $a\in G$, there is an element $a'\in G$ such that $a'\odot a=e$.
\item[3.] For all $a,b\in G$, there is an automorphism ${\rm gyr}[a,b]\in {\rm Aut}(G,\odot)$ such that
$$\hspace{10mm}a\odot(b\odot c)=(a\odot b)\odot {\rm gyr}[a,b](c) \hspace{10mm}\mbox{ (left gyroassociative law)}$$for all $c\in G$.
\item[4.] For all $a,b\in G, {\rm gyr}[a\odot b,b]={\rm gyr}[a,b]. \hspace{10mm}\mbox{ (left loop property)}$
\end{itemize}Moreover, if $$\hspace{37mm}a\odot b={\rm gyr}[a,b](b\odot a) \hspace{10mm}\mbox{ (gyrocommutative law)}$$for all $a,b\in G$,
then $(G,\odot)$ is called a \textit{gyrocommutative gyrogroup}.

We remark that the axioms in the above definition imply the right counterparts. In particular,
any gyrogroup has a unique two-sided identity $e$, and an element $a$ of the gyrogroup
has a unique two-sided inverse $a'$. Given two elements $a,b$ of a gyrogroup $G$, the
map ${\rm gyr}[a,b]$ is called the \textit{gyroautomorphism generated by a and b}. By Theorem 2.10
of \cite{7}, the gyroautomorphisms are completely determined by the gyrator identity
$${\rm gyr}[a,b](c)=(a\odot b)'\odot\left[a\odot(b\odot c)\right]$$for all $a,b,c\in G$. Obviously,
every group forms a gyrogroup under the same operation by defining the gyroautomorphisms to be the identity
automorphism, but the converse is not in general true. From this point of view, gyrogroups suitably generalize groups.

Recall also that \textit{gyrogroup homomorphism} is a map between gyrogroups that preserves the
gyrogroup operations. A bijective gyrogroup homomorphism is called a \textit{gyrogroup
isomorphism}. We say that two gyrogroups $G_1$ and $G_2$ are \textit{isomorphic}, written $G_1\cong G_2$, if
there exists a gyrogroup isomorphism from $G_1$ to $G_2$. Given a gyrogroup $G$, a gyrogroup isomorphism from $G$ to itself
is called a \textit{gyrogroup automorphism} of $G$.

One of the most interesting purely algebraic classes of gyrogroups is introduced in \cite{1}, as follows.
Define on a group $(G,\cdot)$ the binary operation $$a\odot b=a^2ba^{-1}, \forall\, a,b\in G.$$Then,
by Theorem 3.7 of \cite{1}, we have:

\bigskip\noindent{\bf Theorem 1.} {\it $(G,\odot)$ is a gyrogroup if and only if $(G,\cdot)$ is central by a
$2$-Engel group.}
\bigskip

In what follows we will call $(G,\odot)$ the \textit{gyrogroup associated to a given group $(G,\cdot)$, which is assumed to be central by a $2$-Engel group}. Note that in this case
the gyroautomorphism generated by two elements $a$ and $b$ of $G$ is given by $${\rm gyr}[a,b]=\varphi_{[a,b^{-1}]},$$where $\varphi_{[a,b^{-1}]}$
is the inner automorphism of $G$ induced by the commutator $[a,b^{-1}]$ of $a$ and $b^{-1}$.

We are now in a position to characterize the gyrocommutativity of $(G,\odot)$.

\bigskip\noindent{\bf Theorem 2.} {\it $(G,\odot)$ is gyrocommutative if and only if the inner automorphism group of $(G,\cdot)$ is of exponent $3$.}
\bigskip

Clearly, if $(G,\cdot)$ is commutative, then the binary operations $\cdot$ and $\odot$ coincide, and $(G,\odot)$ is gyrocommutative. Note that there is also a non-
commutative group, which is central by a 2-Engel group, such that its associated gyrogroup is gyrocommutative (e.g. the group of upper triangular matrices over $\mathbb{F}_3$ 
with diagonal (1,1,1)).

Next, let $(G_1,\cdot)$ and $(G_2,\cdot)$ be two finite groups central by a 2-Engel group, and let $(G_1,\odot)$ and $(G_2,\odot)$ be their associated gyrogroups. Obviously, if $(G_1,\cdot)\cong(G_2,\cdot)$, then a group isomorphism from $G_1$ to $G_2$ is also a gyrogroup isomorphism from $(G_1,\odot)$ to $(G_2,\odot)$, that is, $(G_1,\odot)\cong(G_2,\odot)$. A sufficient condition for the converse to be true is given in the following theorem.

\bigskip\noindent{\bf Theorem 3.} {\it If $3\nmid |G_1|$, then $(G_1,\cdot)\cong(G_2,\cdot)$ if and only if $(G_1,\odot)\cong(G_2,\odot)$.}
\bigskip

From Theorem 3 we obtain the following corollary.

\bigskip\noindent{\bf Corollary 4.} {\it Let $(G,\cdot)$ be a group central by a 2-Engel group such that $3\nmid |G|$. Then the group of all gyrogroup automorphisms of $(G,\odot)$
coincides with the group of all automorphisms of $(G,\cdot)$.}
\bigskip

Finally, we observe that there is an interesting conservative functor $F$ between the category $\calc$ of finite groups central by a 2-Engel group whose order is not divisible by $3$
and the category of gyrogroups, which associates to each object $(G,\cdot)$ in $\calc$ the gyrogroup $(G,\odot)$ and to each homomorphism $f$ in $\calc$ the gyrogroup homomorphism $F(f)=f$.

Much of our notation is standard and will usually not be repeated
here. Elementary notions and results on groups can be found in
\cite{2}.

\section{Proofs of the main results}

\noindent{\bf Proof of Theorem 2.} In our case the gyrocommutative law becomes $$a^2ba^{-1}=\varphi_{[a,b^{-1}]}(b^2ab^{-1}), \forall\, a,b\in G,$$which means
$$b^3a=ab^3, \forall\, a,b\in G,$$i.e. $$b^3\in Z(G), \forall\, b\in G.$$Obviously, this condition is equivalent to $\exp({\rm Inn}(G))=3$ in view of the group isomorphism $G/Z(G)\cong{\rm Inn}(G)$.
\hfill\rule{1,5mm}{1,5mm}

\bigskip\noindent{\bf Proof of Theorem 3.} Assume that $(G_1,\odot)\cong(G_2,\odot)$ and let $f:G_1\longrightarrow G_2$ be a gyrogroup isomorphism. Then $f$ is a bijection and
$$f(a^2ba^{-1})=f(a)^2f(b)f(a)^{-1}, \forall\, a,b\in G.\0(1)$$If $e_i$ is the identity of $G_i$, $i=1,2$, we infer that $$f(e_1)=e_2$$by taking $a=b=e_1$ in (1). Also, by taking $b=a$ and $b=a^{-1}$ in (1), respectively, one obtains $$f(a^2)=f(a)^2 \mbox{ and } f(a^{-1})=f(a)^{-1}, \forall\, a\in G_1.$$Next, let us write (1) with $a^{-1}ba$ and $a^{-1}b^{-1}$ instead of $b$, respectively. It follows that $$f(ab)=f(a)^2f(a^{-1}ba)f(a)^{-1}\0(2)$$and $$f(aba^{-1})=f(a)f(ba)f(a)^{-2}.\0(3)$$Replace $a$ with $a^{-1}$ in (3). Then $$f(a^{-1}ba)=f(a)^{-1}f(ba^{-1})f(a)^2,$$which together with (2) leads to $$f(ab)=f(a)f(ba^{-1})f(a).$$By writing this equality with $ba$ instead of $b$, we find $$f(aba)=f(a)f(b)f(a), \forall\, a,b\in G_1.\0(4)$$Finally, replacing $b$ with $a^2ba^{-1}$ in (4) gives $$f(a^3b)=f(a)^3f(b), \forall\, a,b\in G_1,\0(5)$$ and taking $b=e_1$ in (5) gives $$f(a^3)=f(a)^3, \forall a\in G_1.\0(6)$$

We are now in a position to prove that $f$ is a group homomorphism. Let $x,y\in G_1$. Since $3\nmid |G_1|$, we have $3\nmid o(x)$ and consequently $\gcd(3,o(x))=1$, i.e. $1=3\alpha+o(x)\beta$ for some integers $\alpha$ and $\beta$. It follows that $$x=x^{3\alpha+o(x)\beta}=x^{3\alpha}x^{o(x)\beta}=x^{3\alpha}.$$Then (5) and (6) lead to $$f(xy)=f(x^{3\alpha}y)=f((x^{\alpha})^3y)=f(x^{\alpha})^3f(y)=f(x^{3\alpha})f(y)=f(x)f(y),$$as desired. Hence $f$ is a group isomorphism, completing the proof.
\hfill\rule{1,5mm}{1,5mm}

\vspace*{5ex}\small

\hfill
\begin{minipage}[t]{5cm}
Marius T\u arn\u auceanu \\
Faculty of  Mathematics \\
``Al.I. Cuza'' University \\
Ia\c si, Romania \\
e-mail: {\tt tarnauc@uaic.ro}
\end{minipage}

\end{document}